\documentclass[11pt,a4paper]{article}
\usepackage{amsfonts}
\usepackage{}

\usepackage{graphicx}
\usepackage{amsfonts,amsmath, amssymb}

\newtheorem{theo}{\textbf{\ \ \quad Theorem}}[section]

\newtheorem{lem}{\textbf{\ \ \quad Lemma}}[section]
\newtheorem{remark}{\textbf{\ \ \quad Remark}}[section]

\newtheorem{prop}{\textbf{\ \ \quad Proposition}}[section]

\newcommand{\lbl}[1]{\label{#1}}

\newcommand{\be}{\begin{equation}}
\newcommand{\ee}{\end{equation}}
\newcommand\bes{\begin{eqnarray}}
\newcommand\ees{\end{eqnarray}}
\newcommand{\bess}{\begin{eqnarray*}}
\newcommand{\eess}{\end{eqnarray*}}

\newcommand{\nm}{\nonumber}

\newcommand{\R}{\mathbb{R}}

{\setlength\arraycolsep{2pt}}

\title{BMO estimates for stochastic singular integral operators and its application to PDEs with L\'{e}vy noise}

\author{Guangying Lv$^a$, Hongjun Gao$^b$ Jinlong Wei$^c$, Jiang-Lun Wu$^d$\\
\\
\ \\
   {\small \it $^a$Institute of Contemporary Mathematics, Henan University}\\
  {\small \it Kaifeng, Henan 475001, China}\\
  {\small \tt gylvmaths@henu.edu.cn}\\
    {\small \it $^b$Institute of Mathematics, School of Mathematical Science}\\
  {\small \it Nanjing Normal University, Nanjing 210023, China}\\
  {\small \tt gaohj@njnu.edu.cn}\\
  {\small \it $^c$ School of Statistics and Mathematics, Zhongnan University of}\\
  {\small \it
 Economics and Law, Wuhan, Hubei 430073, China}\\
   {\small \tt  weijinlong.hust@gmail.com }\\
   {\small \it $^d$ Department of Mathematics, Swansea University, Swansea SA2 8PP, UK}\\
   {\small \tt  j.l.wu@swansea.ac.uk }
}

\begin{document}
\maketitle

\medskip

\begin{abstract}
In this paper, we consider the stochastic singular integral operators and
obtain the BMO estimates. As an application, we consider the fractional Laplacian equation with additive noises
   \bess
du_t(x)=\Delta^{\frac{\alpha}{2}}u_t(x)dt+\sum_{k=1}^\infty\int_{\mathbb{R}^m}g^k(t,x)z\tilde N_k(dz,dt),\ \ \ u_0=0,\ 0\leq t\leq T,
   \eess
where $\Delta^{\frac{\alpha}{2}}=-(-\Delta)^{\frac{\alpha}{2}}$, and $\int_{\mathbb{R}^m}z\tilde N_k(t,dz)=:Y_t^k$ are independent $m$-dimensional pure
jump L\'{e}vy processes with L\'{e}vy measure of $\nu^k$.
Following the idea of \cite{Kim}, we obtain the $q$-th order BMO quasi-norm of the
$\frac{\alpha}{q_0}$-order derivative of $u$ is controlled by the norm of $g$.



\end{abstract}

\baselineskip=15pt

\section{Introduction}
\setcounter{equation}{0}
For a stochastic process $\{X_t,t\in T\}$, there are two important facts worth studying.
One is its probability density function (PDF) or its probability law, the other
is the estimates of moment. But for a stochastic process depending on spatial variable, that
is, $X_t=X(t,\omega,x)$ ($x$ is the spatial variable), it is hard to consider its PDF or
probability law. Fortunately, we can get some estimates of moment. In this paper, we
focus on the estimates of solutions of
stochastic partial differential equations (SPDEs).

For SPDEs,
many kinds of estimates
of the solutions have been well studied. By using parabolic Littlewood-Paley inequality,
 Krylov \cite{Krylov} proved that for
SPDEs of the type
   \bes
du=\Delta udt+gdw_t,
  \lbl{1.1}\ees
it holds that
   \bes
\mathbb{E}\|\nabla u\|_{L^p((0,T)\times\mathbb{R}^d)}^p\leq C(d,p)
\mathbb{E}\|g\|_{L^p((0,T)\times\mathbb{R}^d)}^p,
    \lbl{1.2} \ees
where $w_t$ is a Wiener process and $p\in[2,\infty)$. van Neerven et al. \cite{NVW} introduce
a significant extension of (\ref{1.2}) to a class of operators $A$ which admit a bounded
$H^\infty$-calculus of angle less than $\pi/2$. Kim \cite{Kim} established a BMO estimate
 for stochastic singular integral operators. And as an application, they considered (\ref{1.1}) and obtained
the $q$-th order BMO quasi-norm of the
derivative of $u$ is controlled by $\|g\|_{L^\infty}$. Just recently, Kim et al. \cite{KKL}
studied the parabolic Littlewood-Paley inequality for
a class of time-dependent pseudo-differential operators of arbitrary order,
and applied this result to the high-order stochastic PDE.

Recently, Yang \cite{Y} considered the following SPDEs
 \bess
du=\Delta^{\frac{\alpha}{2}}udt+fdX_t,\ \ \ u_0=0,\ 0< t<T,
  \eess
where $\Delta^{\frac{\alpha}{2}}=-(-\Delta)^{\frac{\alpha}{2}}$, $0<\alpha<2$, and
$X_t$ is a L\'{e}vy process. They obtained a parabolic Triebel-Lizorkin space estimate for the convolution
operator.

Regarding elliptic and parabolic singular integral operators, the BMO estimates was
already established in \cite{Cbook,Gbook}. In this paper, we consider the stochastic singular integral operator
    \bes
\mathcal {G}g(t,x)&=&\int_0^t\int_ZK(t,s,\cdot)\ast g(s,\cdot,z)(x)\tilde N(dz,ds)\nm\\
&=&\int_0^t\int_Z\int_{\mathbb{R}^d}K(t,s,x-y)g(s,y,z)dy\tilde N(dz,ds).
   \lbl{1.3}\ees
Our main purpose is to present appropriate conditions on the kernel $K$ for the following
estimate:
    \bes
[\mathcal {G}g]_{\mathbb{BMO}(T,q)}&\leq&
N\left(\left\|\left(\int_Z
\|g(\cdot,\cdot,z)\|^2_{L^\infty(\mathcal {O}_T)}\nu(dz)\right)^{q/2}\right\|_{L^{\frac{\kappa}{\kappa-q}}}
\right.\nm\\
&&\qquad+\left\|\int_Z|g(\cdot,\cdot,z)|^{q_0}_{L^\infty(\mathcal {O}_T)}\nu(dz)\right\|_{L^{\tilde\kappa}(\Omega)}^{q/q_0}\nm\\
&&\qquad\left.
+\left\|\int_Z\|g(\cdot,\cdot,z)\|^q_{L^\infty(\mathcal {O}_T)}\nu(dz)\right\|_{L^{\frac{\kappa}{\kappa-q}}}\right),
  \lbl{1.4}\ees
where $q\in[2,p_0\wedge\kappa]$, $\tilde\kappa$ is the conjugate of a positive constant $\kappa$, the constant $N$ depends on $q$ and $d$,
and $\nu$ is a measure, see Section 2. As an application of (\ref{1.4}), we prove that the solution of
the following equation
   \bess
du_t(x)=\Delta^{\frac{\alpha}{2}}u_t(x)dt+\sum_{k=1}^\infty\int_{\mathbb{R}^m}g^k(t,x)z\tilde N_k(dz,dt),\ \ \ u_0=0,\ 0\leq t\leq T,
   \eess
satisfies that for $q\in[2,q_0]$
      \bess
[\nabla^{\beta}u]_{\mathbb{BMO}(T,q)}&\leq& N\hat c\left(\mathbb{E}[\||g|_{\ell_2}\|^{q_0}_{L^\infty(\mathcal {O}_T)}]\right)^{q/q_0},
  \eess
where $\int_{\mathbb{R}^m}z\tilde N_k(t,dz)=:Y_t^k$ are independent $m$-dimensional pure
jump L\'{e}vy processes with L\'{e}vy measure of $\nu^k$, $\beta=\alpha/q_0$ and $\hat c$ is defined as in (\ref{4.3}),
see Section 4 for details. Moreover, we find if we consider the following stochastic parabolic equation
   \bess
du_t(x)=\Delta^{\frac{\alpha}{2}}u_t(x)dt+\sum_{k=1}^\infty h^k(t,x)dW_t^k,\ \ \ u_0=0,\ 0\leq t\leq T,
   \eess
where $W_t^k$ are independent one-dimensional Wiener processes. We have the following estimate, for any $q\in(0,p]$,
      \bess
[\nabla^{\frac{\alpha}{2}}u]_{\mathbb{BMO}(T,q)}&\leq& N\left(\mathbb{E}[\||h|_{\ell_2}\|^{p}_{L^\infty(\mathcal {O}_T)}]\right)^{1/p}.
  \eess
under the condition that $h\in L^p(T,\ell_2)$, see Theorem \ref{t4.2}. Specially, taking $\alpha=2$, we obtain
the result of \cite[Theorem 3.4]{Kim}.

Due to the difference between the Brownian motion and L\'{e}vy process, it is more difficult
to get the BMO estimate for L\'{e}vy process. Following the idea of \cite{Kim}, we obtain
the BMO estimate of stochastic singular integral operators. We remark that there are many
places different from those in \cite{Kim}. First, the assumptions on the kernel are different from
those in \cite{Kim}, see Section 2; Second, the exponent $q$ in \cite{Kim} do not depend on the
properties of kernel but we do. For simplicity, we only consider a simple case, see the discussion in Section 4.

This paper is organized as follows. In Section 2, we introduce the main results. The proof of the main
results is complete in section 3. Section 4 is concerned with an application of our result. This paper ends
with a short discussion, which shows that we can give a simple proof of the result in Section 2 if
the function $g$ has high regularity.

Before we end this section, we introduce some notations used in this paper.
As usual $\mathbb{R}^d$ stands for the Euclidean space of points $x=(x_1,\cdots,x_d)$,
$B_r(x):=\{y\in\mathbb{R}^d:|x-y|<r\}$ and $B_r:=B_r(0)$. $\mathbb{R}_+$ denotes the
set $\{x\in\mathbb{R},x>0\}$. $a\wedge b=\min\{a,b\}$, $a\vee b=\max\{a,b\}$ and
$L^p:=L^p(\mathbb{R}^d)$.
$N=N(a,b,\cdots)$ means that the constant
$N$ depends only on $a,b,\cdots$.

\section{Known results and Main result}
\setcounter{equation}{0}

Let $(\Omega,\mathcal {F},\mathbb{F},\mathbb{P})$ be a complete probability space such that
$\mathcal {F}_t$ is a filtration on $\Omega$ containing
all $P$-null subsets of $\Omega$ and $\mathbb{F}$ be the predictable $\sigma$-field by
$(\mathcal {F}_t,t\geq0)$. We are given a measure space $(Z,\mathcal {Z},\nu)$ and a Poisson
measure $\mu$ on $[0,T]\times Z$, defined on the stochastic basis. The compensator of
$\mu$ is Leb$\otimes\nu$, and the compensated measure $\tilde N:=\mu-Leb\otimes\nu$

Fix $\gamma>0$ and $T\in(0,\infty]$. Denote
   \bess
\mathcal {O}_T=(0,T)\times\mathbb{R}^d.
  \eess
For a measurable function $h$ on $\Omega\times\mathcal {O}_T$, we define the $q$-th
order stochastic BMO (Bounded mean oscillation) quasi-norm of $h$ on $\Omega\times\mathcal {O}_T$
as follows:
   \bess
[h]_{\mathbb{BMO}(T,q)}^q=\sup_Q\frac{1}{Q^2}\mathbb{E}\int_Q\int_Q|h(t,x)-h(s,y)|^qdtdxdsdy,
   \eess
where the sup is taken over all $Q$ of the type
   \bess
Q=Q_c(t_0,x_0):=(t_0-c^\gamma,t_0+c^\gamma)\times B_c(x_0)\subset\mathcal {O}_T,\ \ c>0,t_0>0.
  \eess
It is remarked that when $q=1$, this is equivalent to the classical BMO semi-norm which is
introduced by John-Nirenberg \cite{JN}.

Let $K(\omega,t,s,x)$ be a measurable function on $\Omega\times\mathbb{R}_+\times\mathbb{R}_+\times\mathbb{R}^d$
such that for each $t\in\mathbb{R}_+$, $(\omega,s)\mapsto K(\omega,t,s,\cdot)$ is a predictable
$L_{loc}^1$-valued process.

Firstly, we recall the results of \cite{Kim}. In \cite{Kim}, the following assumptions are needed.

{\bf Assumption 2.1} There exist a $\kappa\in[1,\infty]$ and a nondecreasing
function $\varphi(t):(0,\infty)\mapsto[0,\infty)$ such that

(i) for any $t>\lambda>0$ and $c>0$,
   \bess
\left\|\int_\lambda^t\Big|\int_{|x|\geq c}|K(t,r,x)|dx\Big|^2dr\right\|_{L^{\kappa/2}(\Omega)}\leq\varphi((t-\lambda)c^{-\gamma});
   \eess

(ii) for any $t>s>\lambda>0$,
   \bess
\left\|\int_0^\lambda\left(\int_{\mathbb{R}^d}|K(t,r,x)-K(s,r,x)|dx\right)^2dr\right\|_{L^{\kappa/2}(\Omega)}
\leq\varphi((t-s)(t\wedge s-\lambda)^{-1});
   \eess

(iii) for any $s>\lambda\geq0$ and $h\in\mathbb{R}^d$,
    \bess
\left\|\int_0^\lambda\left(\int_{\mathbb{R}^d}|K(s,r,x+h)-K(s,r,x)|dx\right)^2dr\right\|_{L^{\kappa/2}(\Omega)}
\leq N\varphi(|h|(s-\lambda)^{-1/\gamma}).
   \eess

{\bf Assumption 2.2} Suppose that $\mathcal {G}g$ is well-defined (a.e.) and the following holds:
  \bess
\mathbb{E}\int_0^T\|\mathcal {G}g(t,\cdot)\|_{L^{p_0}}^{p_0}dt
\leq N_0\left\|\int_0^T\||g(t,\cdot)|_{l_2}\|_{L^{p_0}}^{p_0}dt\right\|_{L^{\tilde\kappa}(\Omega)},
   \eess
where $\tilde\kappa$ is the conjugate of $\kappa$, and
 \bess
\mathcal {G}g(t,x)=\sum_{k=1}^\infty\int_0^t\int_{\mathbb{R}^d}K(t,s,x-y)g^k(s,y)dydw_s^k,
   \eess
with $w_t$ is a Wiener process.
Under the Assumptions 2.1 and 2.2, Kim obtained
the BMO estimate of $\mathcal {G}g$.

Comparing with the assumption 2.1, due to the Kunita's first inequality (see Page 265 of \cite{Ap}),
we need the following assumptions. For the Kunita's inequality and BDG inequality of L\'{e}vy noise,
see Lemma 3.1 of \cite{MPR} and \cite{MR} respectively.

{\bf Assumption 2.3} There exist constants $q_0\geq2$,
 $\kappa\in[1,\infty]$ and a nondecreasing
function $\varphi(t):(0,\infty)\mapsto[0,\infty)$ such that

(i) for any $t>\lambda>0$ and $c>0$,
   \bess
\left\|\int_\lambda^t\Big|\int_{|x|\geq c}|K(t,r,x)|dx\Big|^{q_0}dr\right\|_{L^{\kappa/q_0}(\Omega)}\leq\varphi((t-\lambda)c^{-\gamma});
   \eess

(ii) for any $t>s>\lambda>0$,
   \bess
\left\|\int_0^\lambda\left(\int_{\mathbb{R}^d}|K(t,r,x)-K(s,r,x)|dx\right)^{q_0}dr\right\|_{L^{\kappa/q_0}(\Omega)}
\leq\varphi((t-s)(t\wedge s-\lambda)^{-1});
   \eess

(iii) for any $s>\lambda\geq0$ and $h\in\mathbb{R}^d$,
    \bess
\left\|\int_0^\lambda\left(\int_{\mathbb{R}^d}|K(s,r,x+h)-K(s,r,x)|dx\right)^{q_0}dr\right\|_{L^{\kappa/q_0}(\Omega)}
\leq N\varphi(|h|(s-\lambda)^{-1/\gamma}).
   \eess

\begin{remark}\lbl{r2.1} The difference between assumptions 2.1 and 2.3 is because
the following Kunita's first inequality.
   \bes
\mathbb{E}\left(\sup_{0\leq t\leq T}|I(t)|^p\right)&\leq&N(p)\left\{
\mathbb{E}\left[\left(\int_0^T\int_Z|H(t,z)|^2\nu(dz)dt\right)^{p/2}\right]\right.\nm\\
&&\left.+\mathbb{E}\left[\int_0^T\int_Z|H(t,z)|^p\nu(dz)dt\right]\right\},
   \lbl{2.1}\ees
where $p\geq2$ and
    \bess
I(t)=\int_0^t\int_ZH(s,z)\tilde N(dz,ds).
   \eess
When $\tilde N(ds,dz)$ is replaced by $dw_sdz$, the second term of right
side hand of (\ref{2.1}) will disappear. Hence, in order to deal with
the difficult from the L\'{e}vy process, we give the assumption 2.3.
   \end{remark}

 {\bf Assumption 2.4} Similar  Assumption 2.2, suppose that $\mathcal {G}g$ is well-defined (a.e.)
 and the following holds:
   \bes
\mathbb{E}\int_0^T\|\mathcal {G}g(t,\cdot)\|_{L^{q_0}}^{q_0}dt
\leq N_0\left\|\int_0^T\int_Z\|g(t,\cdot,z)\|_{L^{q_0}}^{q_0}\nu(dz)dt\right\|_{L^{\tilde\kappa}(\Omega)}.
   \lbl{2.2}\ees

Our main result is the following.
\begin{theo}\lbl{t2.1} Let Assumptions 2.3 and 2.4 hold. Assume that the function
$g$ satisfies
   \bes
 \left\|\int_Z
\|g(\cdot,\cdot,z)\|^\varpi_{L^\infty(\mathcal {O}_T)}\nu(dz)\right\|_{L^{\varsigma}(\Omega)}<\infty,\ \
\ \varpi=2\ \ {\rm or}\ \ q_0,
  \lbl{j.20}\ees
where $\varsigma=q_0\tilde\kappa\vee\frac{q_0\kappa}{2(\kappa-q_0)^+}$ ($\varsigma=\infty$ if $\kappa\leq q_0$).
Then for any $q\in[2,q_0\wedge \kappa]$,
one has
     \bes
[\mathcal {G}g]_{\mathbb{BMO}(T,q)}&\leq&
N\left(\left\|\left(\int_Z
\|g(\cdot,\cdot,z)\|^2_{L^\infty(\mathcal {O}_T)}\nu(dz)\right)^{q/2}\right\|_{L^{\frac{\kappa}{\kappa-q}}}
\right.\nm\\
&&\qquad+\left\|\int_Z\|g(\cdot,\cdot,z)\|^{q_0}_{L^\infty(\mathcal {O}_T)}\nu(dz)\right\|_{L^{\tilde\kappa}(\Omega)}^{q/q_0}\nm\\
&&\qquad\left.
+\left\|\int_Z\|g(\cdot,\cdot,z)\|^q_{L^\infty(\mathcal {O}_T)}\nu(dz)\right\|_{L^{\frac{\kappa}{\kappa-q}}}\right),
  \lbl{2.3}\ees
where $N=N(N_0,d,q,q_0,\gamma,\kappa,\varphi)$.
\end{theo}

\begin{remark}\lbl{r2.2} 1. Comparing Theorem \ref{t2.1} with Theorem 2.4 in \cite{Kim}, it
is not hard to find in Theorem 2.4 of \cite{Kim} the exponent $q$ does not depend on $q_0$.
Actually, the range of exponent $q$ is $(0,p_0\wedge \kappa]$ and in this paper is $[2,q_0\wedge \kappa]$.
In other words, the range of exponent $q$ depends on the properties of kernel $K$.
The lower bound of $q$ is because the Kunita's first inequality holds for $q\geq2$.

2. In Theorem \ref{t2.1}, we did not write the right hand of (\ref{2.3}) as
a uniform format. The reason is that $\int_Z\nu(dz)$ maybe not exist. If we assume that
   \bess
\int_Z (z^2\wedge1)\nu(dz)\leq N_1\,\ {\rm and}\ \
\left\|\int_Z\|g(\cdot,\cdot,z)\|^{q_0}_{L^\infty(\mathcal {O}_T)}(1+f(z)^{-\frac{q_0}{2}}+f(z)^{-\frac{q_0}{q}})\nu(dz)\right\|_{L^{\kappa^*}(\Omega)}<\infty,
    \eess
where $N_1$ is a positive constant, then (\ref{2.3}) can be replaced by
   \bess
[\mathcal {G}g]_{\mathbb{BMO}(T,q)}\leq
\left\|\int_Z\|g(\cdot,\cdot,z)\|^{q_0}_{L^\infty(\mathcal {O}_T)}(1+f(z)^{-\frac{q_0}{2}}+f(z)^{-\frac{q_0}{q}})\nu(dz)\right\|_{L^{\kappa^*}(\Omega)}^{q/q_0},
   \eess
where
   \bess
\kappa^*=\tilde\kappa\vee\frac{\kappa}{\kappa-q},\ \ f(z)=\frac{z^2+1-|z^2-1|}{2}=z^2\wedge1.
   \eess

3. The condition (\ref{j.20}) coincides with (\ref{4.3}) in Section 4. Under the condition (\ref{j.20}),
it is easy to check that
   \bess
 \left\|\int_Z
\|g(\cdot,\cdot,z)\|^q_{L^\infty(\mathcal {O}_T)}\nu(dz)\right\|_{L^{\frac{\kappa}{\kappa-q}}(\Omega)}<\infty.
   \eess
\end{remark}

\section{Proof of the main result}
\setcounter{equation}{0}

In this section, we first estimate the expectation of local mean average of
$\mathcal {G}g$ and its difference in terms of the supremum of $|g|$ given a
vanishing condition on $g$. Then we complete the proof of main result.

\begin{lem}\lbl{l3.1} Let $q\in[2,q_0]$, $0\leq a\leq b\leq T$, and Assumption 2.4 hold.
Suppose that $g$ vanishes on $(a,b)\times(B_{3c})^c\times Z$ and $(0,a)\times\mathbb{R}^d\times Z$. Then
   \bess
\mathbb{E}\int_a^b\int_{B_c}|\mathcal {G}g(t,x)|^qdxdt\leq
N(b-a)|B_{3c}|\left\|\sup_{(a,b)\times B_{3c}}\int_Z|g(\cdot,\cdot,z)|^{q_0}\nu(dz)\right\|_{L^{\tilde\kappa}(\Omega)}^{q/p_0},
  \eess
where $N=N(N_0)$.
  \end{lem}

{\bf Proof.} The proof of this lemma is similar to that of Lemma 4.1
in \cite{Kim}. In order to read easily, we give the outline of the proof.
By H$\ddot{o}$lder's inequality and Assumption 2.4,
   \bess
&&\mathbb{E}\int_a^b\int_{B_c}|\mathcal {G}g(t,x)|^qdxdt\\
&\leq&(b-a)^{(q_0-q)/q_0}|B_c|^{(q_0-q)/q_0}\left(\mathbb{E}\int_a^b\int_{B_c}|\mathcal {G}g(t,x)|^{q_0}dxdt\right)^{q/q_0}\\
&\leq&N(b-a)^{(q_0-q)/q_0}|B_c|^{(q_0-q)/q_0}
\left\|\int_0^T\int_Z\|g(t,\cdot,z)\|_{L^{q_0}}^{q_0}\nu(dz)dt\right\|_{L^{\tilde\kappa}(\Omega)}^{q/q_0}.
   \eess
Since $g$ vanishes on $(a,b)\times(B_{3c})^c$ and $(0,a)\times\mathbb{R}^d$, the above term is
equal to or less than
   \bess
&&N(b-a)^{(q_0-q)/q_0}|B_c|^{(q_0-q)/q_0}
\left\|\int_a^b\int_{B_{3c}}\int_Z|g(t,x,z)|^{q_0}\nu(dz)dxdt\right
\|_{L^{\tilde\kappa}(\Omega)}^{q/q_0}\\
&\leq&N(b-a)|B_{3c}|\left\|\sup_{(a,b)\times B_{3c}}\int_Z|g(\cdot,\cdot,z)|^{q_0}\nu(dz)\right\|_{L^{\tilde\kappa}(\Omega)}^{q/q_0}.
  \eess
The proof of lemma is complete. $\Box$

\begin{lem}\lbl{l3.2} Let $q\in[2,q_0\wedge \kappa]$, $0\leq a\leq b\leq T$ and Assumption 2.3 (i) hold.
Suppose that $g$ vanishes on $(0,\frac{3b-a}{2})\times B_{2c}\times Z$. Then
   \bes
&&\mathbb{E}\int_a^b\int_{B_c}\int_a^b\int_{B_c}|\mathcal {G}g(t,x)-\mathcal {G}g(s,y)|^qdxdtdsdy\nm\\
&\leq&N(b-a)^2|B_{c}|^2[\varphi(bc^{-\gamma})]^{q/q_0}\left(\left\|\left(\int_Z
\|g(\cdot,\cdot,z)\|^2_{L^\infty(\mathcal {O}_T)}\nu(dz)\right)^{q/2}\right\|_{L^{\frac{\kappa}{\kappa-q}}}
\right.\nm\\
&&\left.
+\left\|\int_Z\|g(\cdot,\cdot,z)\|^q_{L^\infty(\mathcal {O}_T)}\nu(dz)\right\|_{L^{\frac{\kappa}{\kappa-q}}}\right),
  \lbl{3.1}\ees
where $\frac{\infty}{\infty}:=1$ and $N=N(T,q)$.
  \end{lem}

{\bf Proof.} Let $(t,x)\in(a,b)\times B_c$ and $0\leq r\leq t$. If $|y|\leq c$, then
$(r,x-y)\in(0,\frac{3b-a}{2})\times B_{2c}$ and $g(r,x-y,z)=0$ for all $z\in Z$. Hence,
Assumption 2.3 (i), H$\ddot{o}$lder inequality and Kunita's first inequality (\ref{2.1}) implies
   \bess
\mathbb{E}|\mathcal {G}g(t,x)|^q&\leq&
\mathbb{E}\left(\int_0^t\int_Z|\int_{\mathbb{R}^d}K(t,r,y)g(r,x-y,z)dy|^2\nu(dz)dr\right)^{q/2}\\
&&+\mathbb{E}\left(\int_0^t\int_Z|\int_{\mathbb{R}^d}K(t,r,y)g(r,x-y,z)dy|^q\nu(dz)dr\right)\\
&\leq&
\mathbb{E}\left(\int_0^t\int_Z|\int_{|y|\geq c}K(t,r,y)g(r,x-y,z)dy|^2\nu(dz)dr\right)^{q/2}\\
&&+\mathbb{E}\left(\int_0^t\int_Z|\int_{|y|\geq c}K(t,r,y)g(r,x-y,z)dy|^q\nu(dz)dr\right)\\
&\leq&T^{(q_0-2)q/(2q_0)}\mathbb{E}\left[\left(\int_0^t\Big|\int_{|y|\geq c}|K(t,r,y)|dy\Big|^{q_0}dr\right)^{q/q_0}\right.\\
&&\qquad\qquad\quad\left.\times
\left(\int_Z\|g(\cdot,\cdot,z)\|^2_{L^\infty(\mathcal {O}_T)}\nu(dz)\right)^{q/2}\right]\\
&&+\mathbb{E}\left[\left(\int_0^t\Big|\int_{|y|\geq c}|K(t,r,y)|dy\Big|^qdr\right)\int_Z\|g(\cdot,\cdot,z)\|^q_{L^\infty(\mathcal {O}_T)}\nu(dz)\right]\\
&\leq&N(T)\left\|\int_0^t\Big|\int_{|y|\geq c}|K(t,r,y)|dy\Big|^{q_0}dr\right\|_{L^{\kappa/q_0}}^{q/q_0}\\
&&\quad\times
\left\|\left(\int_Z\|g(\cdot,\cdot,z)\|^2_{L^\infty(\mathcal {O}_T)}\nu(dz)\right)^{q/2}\right\|_{L^{\frac{\kappa}{\kappa-q}}}\\
&&+N(T)\left\|\int_0^t\Big|\int_{|y|\geq c}|K(t,r,y)|dy\Big|^{q_0}dr\right\|_{L^{\kappa/q_0}}^{q/q_0}\\
&&\quad\times
\left\|\int_Z\|g(\cdot,\cdot,z)\|^q_{L^\infty(\mathcal {O}_T)}\nu(dz)\right\|_{L^{\frac{\kappa}{\kappa-q}}}\\
&\leq&N(T)[\varphi(bc^{-\gamma})]^{q/q_0}\left(\left\|\left(\int_Z\|g(\cdot,\cdot,z)\|^2_{L^\infty(\mathcal {O}_T)}\nu(dz)\right)^{q/2}\right\|_{L^{\frac{\kappa}{\kappa-q}}}
\right.\\
&&\qquad\qquad\qquad\qquad\quad\left.
+\left\|\int_Z\|g(\cdot,\cdot,z)\|^q_{L^\infty(\mathcal {O}_T)}\nu(dz)\right\|_{L^{\frac{\kappa}{\kappa-q}}}\right),
  \eess
which implies that
   \bess
&&\mathbb{E}\int_a^b\int_{B_c}\int_a^b\int_{B_c}|\mathcal {G}g(t,x)-\mathcal {G}g(s,y)|^qdxdtdsdy\nm\\
&\leq&N(q)(b-a)|B_c|\mathbb{E}\int_a^b\int_{B_c}|\mathcal {G}g(t,x)|^qdxdt\\
&\leq&N(T,q)(b-a)^2|B_{c}|^2[\varphi(bc^{-\gamma})]^{q/q_0}\left(\left\|\left(\int_Z\|g(\cdot,\cdot,z)\|^2_{L^\infty(\mathcal {O}_T)}\nu(dz)\right)^{q/2}\right\|_{L^{\frac{\kappa}{\kappa-q}}}
\right.\\
&&\left.
+\left\|\int_Z\|g(\cdot,\cdot,z)\|^q_{L^\infty(\mathcal {O}_T)}\nu(dz)\right\|_{L^{\frac{\kappa}{\kappa-q}}}\right).
  \eess
The inequality (\ref{3.1}) is obtained. The proof of lemma is complete. $\Box$

\begin{lem}\lbl{l3.3} Let $q\in[2,q_0\wedge \kappa]$, $0\leq a< b\leq T$ such that $3a>b$. Suppose that Assumption 2.3 holds
and $\mathcal {G}g$ is well-defined almost everywhere. Assume further that $g$ vanishes on $(\frac{3a-b}{2},\frac{3b-a}{2})\times B_{2c}\times Z$. Then
   \bes
&&\mathbb{E}\int_a^b\int_{B_c}\int_a^b\int_{B_c}|\mathcal {G}g(t,x)-\mathcal {G}g(s,y)|^qdxdtdsdy\nm\\
&\leq&N(b-a)^2|B_{c}|^2\Phi(a,b,c)\left(\left\|\left(\int_Z\|g(\cdot,\cdot,z)\|^2_{L^\infty(\mathcal {O}_T)}\nu(dz)\right)^{q/2}\right\|_{L^{\frac{\kappa}{\kappa-q}}}
\right.\nm\\
&&\left.
+\left\|\int_Z\|g(\cdot,\cdot,z)\|^q_{L^\infty(\mathcal {O}_T)}\nu(dz)\right\|_{L^{\frac{\kappa}{\kappa-q}}}\right),
  \lbl{3.2}\ees
where $N=N(T,q,a,b,c)$ and
   \bess
\Phi(a,b,c)=[\varphi(2)]^{q/q_0}+[\varphi((b-a)c^{-\gamma})]^{q/q_0}+[\varphi(2^{1+1/\gamma}c(b-a)^{-1/\gamma})]^{q/q_0}.
  \eess
  \end{lem}

{\bf Proof.} Due to the Fubini's Theorem, it suffices to prove that for all $(t,x)\in(a,b)\times B_c$ and
 $(s,y)\in(a,b)\times B_c$, the following inequality holds:
    \bess
\mathbb{E}|\mathcal {G}g(t,x)-\mathcal {G}g(s,y)|^q
&\leq&N\Phi(a,b,c)\left(\left\|\left(\int_Z\|g(\cdot,\cdot,z)\|^2_{L^\infty(\mathcal {O}_T)}\nu(dz)\right)^{q/2}\right\|_{L^{\frac{\kappa}{\kappa-q}}}
\right.\nm\\
&&\left.
+\left\|\int_Z\|g(\cdot,\cdot,z)\|^q_{L^\infty(\mathcal {O}_T)}\nu(dz)\right\|_{L^{\frac{\kappa}{\kappa-q}}}\right),
   \eess
Obviously,
   \bess
&&\mathbb{E}|\mathcal {G}g(t,x)-\mathcal {G}g(s,y)|^q\\
&\leq&N(\mathbb{E}|\mathcal {G}g(t,x)-\mathcal {G}g(s,x)|^q+\mathbb{E}|\mathcal {G}g(s,x)-\mathcal {G}g(s,y)|^q)\\
&=:&N(I_1+I_2).
  \eess
{\bf Estimate of $I_1$.} Without loss of generality we assume $t\geq s$.
Hence by Lemma 3.1 of \cite{MPR} and (\ref{2.1}), we get
   \bess
I_1&=&\mathbb{E}|\mathcal {G}g(t,x)-\mathcal {G}g(s,x)|^q\\
&=&\mathbb{E}\left[\Big|\int_0^t\int_Z\int_{\mathbb{R}^d}K(t,r,x-y)g(r,y,z)dy\tilde N(dz,dr)\right.\\
&&\qquad\left.
-\int_0^s\int_Z\int_{\mathbb{R}^d}K(s,r,x-y)g(r,y,z)dy\tilde N(dz,dr)\Big|^q\right]\\
&\leq&N\mathbb{E}\left[\Big|\int_0^t\int_Z\int_{\mathbb{R}^d}K(t,r,x-y)g(r,y,z)dy\tilde N(dz,dr)\right.\\
&&\qquad\left.
-\int_0^s\int_Z\int_{\mathbb{R}^d}K(t,r,x-y)g(r,y,z)dy\tilde N(dz,dr)\Big|^q\right]\\
&&+N\mathbb{E}\left[\Big|\int_0^s\int_Z\int_{\mathbb{R}^d}(K(t,r,x-y)-K(s,r,x-y))g(r,y,z)dy\tilde N(dz,dr)\Big|^q\right]\\
&\leq&N\mathbb{E}\left[\left(\int_s^t\int_Z|\int_{\mathbb{R}^d}K(t,r,x-y)g(r,y,z)dy|^2\nu(dz)dr\right)^{q/2}\right]\nm\\
&& +N\mathbb{E}\left[\int_s^t\int_Z|\int_{\mathbb{R}^d}K(t,r,x-y)g(r,y,z)dy|^q\nu(dz)dr\right] \\
&&+N\mathbb{E}\left[\left(\int_0^s\int_Z|\int_{\mathbb{R}^d}(K(t,r,x-y)-K(s,r,x-y))g(r,y,z)dy|^2\nu(dz)dr\right)^{q/2}\right]\nm\\
&& +N\mathbb{E}\left[\int_0^s\int_Z|\int_{\mathbb{R}^d}(K(t,r,x-y)-K(s,r,x-y))g(r,y,z)dy|^q\nu(dz)dr\right]\\
&=:&N(I_{11}+I_{12}+I_{13}+I_{14}).
  \eess
Note that $g$ vanishes on $(\frac{3a-b}{2},\frac{3b-a}{2})\times B_{2c}\times Z$ and $a>\frac{3a-b}{2}$. {\bf Assumption
2.3 (i)} with $\lambda=s$ yields that
  \bess
I_{11}+I_{12}&=&\mathbb{E}\left[\left(\int_s^t\int_Z|\int_{\mathbb{R}^d}K(t,r,y)g(r,x-y,z)dy|^2\nu(dz)dr\right)^{q/2}\right]\nm\\
&& +\mathbb{E}\left[\int_s^t\int_Z|\int_{\mathbb{R}^d}K(t,r,y)g(r,x-y,z)dy|^q\nu(dz)dr\right] \\
&\leq&\mathbb{E}\left[\left(\int_s^t\Big|\int_{|y|\geq c}|K(t,r,y)|dy\Big|^2dr\int_Z\|g(\cdot,\cdot,z)\|^2_{L^\infty(\mathcal {O}_T)}\nu(dz)\right)^{q/2}\right]\nm\\
&& +\mathbb{E}\left[\int_s^t\Big|\int_{|y|\geq c}|K(t,r,y)|dy\Big|^q\int_Z\|g(\cdot,\cdot,z)\|^q_{L^\infty(\mathcal {O}_T)}\nu(dz)dt\right] \\
&\leq&N[\varphi((b-a)c^{-\gamma})]^{q/q_0}\left(\left\|\left(\int_Z\|g(\cdot,\cdot,z)\|^2_{L^\infty(\mathcal {O}_T)}\nu(dz)\right)^{q/2}\right\|_{L^{\frac{\kappa}{\kappa-q}}}
\right.\\
&&\left.
+\left\|\int_Z\|g(\cdot,\cdot,z)\|^q_{L^\infty(\mathcal {O}_T)}\nu(dz)\right\|_{L^{\frac{\kappa}{\kappa-q}}}\right).
  \eess
Similarly, due to $g$ vanishes on $(\frac{3a-b}{2},\frac{3b-a}{2})\times B_{2c}\times Z$, we divide $(0,s)$ into two parts
$(0,\frac{3a-b}{2})$ and $(\frac{3a-b}{2},s)$. And thus we have
    \bess
I_{13}+I_{14}&=&\mathbb{E}\left[\left(\int_{\frac{3a-b}{2}}^s\int_Z|\int_{\mathbb{R}^d}(K(t,r,x-y)-K(s,r,x-y))g(r,y,z)dy|^2\nu(dz)dr\right)^{q/2}\right]\nm\\
&& +\mathbb{E}\left[\int_{\frac{3a-b}{2}}^s\int_Z|\int_{\mathbb{R}^d}(K(t,r,x-y)-K(s,r,x-y))g(r,y,z)dy|^q\nu(dz)dr\right]\\
&&+\mathbb{E}\left[\left(\int_0^{\frac{3a-b}{2}}\int_Z|\int_{\mathbb{R}^d}(K(t,r,x-y)-K(s,r,x-y))g(r,y,z)dy|^2\nu(dz)dr\right)^{q/2}\right]\nm\\
&& +\mathbb{E}\left[\int_0^{\frac{3a-b}{2}}\int_Z|\int_{\mathbb{R}^d}(K(t,r,x-y)-K(s,r,x-y))g(r,y,z)dy|^q\nu(dz)dr\right] \\
 &=:&I_{131}+I_{141}+I_{132}+I_{142}.
   \eess
Using again {\bf Assumption 2.3 (i)} with $\lambda=\frac{3a-b}{2}$, we get
   \bess
I_{131}+I_{141}&\leq&\mathbb{E}\left[\left(\int_{\frac{3a-b}{2}}^t\int_Z|\int_{\mathbb{R}^d}|K(t,r,x-y)g(r,y,z)|dy|^2\nu(dz)dr\right)^{q/2}\right]\nm\\
&&+\mathbb{E}\left[\left(\int_{\frac{3a-b}{2}}^s\int_Z|\int_{\mathbb{R}^d}|K(s,r,x-y)g(r,y,z)|dy|^2\nu(dz)dr\right)^{q/2}\right]\nm\\
&& +\mathbb{E}\left[\int_{\frac{3a-b}{2}}^t\int_Z|\int_{\mathbb{R}^d}|K(t,r,x-y) g(r,y,z)|dy|^q\nu(dz)dr\right]\\
&& +\mathbb{E}\left[\int_{\frac{3a-b}{2}}^s\int_Z|\int_{\mathbb{R}^d}|K(s,r,x-y) g(r,y,z)|dy|^q\nu(dz)dr\right]\\
&\leq&N[\varphi(2(b-a)c^{-\gamma})]^{q/q_0}\left(\left\|\left(\int_Z\|g(\cdot,\cdot,z)\|^2_{L^\infty(\mathcal {O}_T)}\nu(dz)\right)^{q/2}\right\|_{L^{\frac{\kappa}{\kappa-q}}}
\right.\\
&&\left.
+\left\|\int_Z\|g(\cdot,\cdot,z)\|^q_{L^\infty(\mathcal {O}_T)}\nu(dz)\right\|_{L^{\frac{\kappa}{\kappa-q}}}\right).
 \eess
On the other hand, {\bf Assumption 2.3 (ii)} with $\lambda=\frac{3a-b}{2}$ gives
   \bess
I_{132}+I_{142}&\leq&
N\mathbb{E}\left[\left(\int_0^{\frac{3a-b}{2}}\Big|\int_{\mathbb{R}^d}|K(t,r,x-y)-K(s,r,x-y)|dy\Big|^2dr\right.\right.\\
&&\left.\left.\times
\int_Z\|g(\cdot,\cdot,z)\|^2_{L^\infty(\mathcal {O}_T)}\nu(dz)\right)^{q/2}\right]\nm\\
&& +\mathbb{E}\left[\int_0^{\frac{3a-b}{2}}\Big|\int_{\mathbb{R}^d}|K(t,r,x-y)-K(s,r,x-y)|dy\Big|^q\right.\\
&&\left.\times
\int_Z\|g(\cdot,\cdot,z)\|^q_{L^\infty(\mathcal {O}_T)}\nu(dz)dr\right] \\
&\leq&N[\varphi(2)]^{q/q_0}\left(\left\|\left(\int_Z\|g(\cdot,\cdot,z)\|^2_{L^\infty(\mathcal {O}_T)}\nu(dz)\right)^{q/2}\right\|_{L^{\frac{\kappa}{\kappa-q}}}
\right.\\
&&\left.
+\left\|\int_Z\|g(\cdot,\cdot,z)\|^q_{L^\infty(\mathcal {O}_T)}\nu(dz)\right\|_{L^{\frac{\kappa}{\kappa-q}}}\right),
  \eess
where we used $s-\frac{3a-b}{2}\geq a-\frac{3a-b}{2}=\frac{b-a}{2}$ and $(t-s)(s-\frac{3a-b}{2})^{-1}\leq2$.

{\bf Estimate of $I_2$.} By using the fact $g=0$ on $(\frac{3a-b}{2},\frac{3b-a}{2})\times B_{2c}\times Z$ again, we divide $(0,s)$ into two parts
$(0,\frac{3a-b}{2})$ and $(\frac{3a-b}{2},s)$.
Direct calculations shows that
   \bess
I_2&\leq&N \mathbb{E}\left(\int_0^s\int_Z\Big|\int_{\mathbb{R}^d}K(s,r,w)(g(r,x-w,z)-g(r,y-w,z))dw|^2\nu(dz)dr\right)^{q/2}\\
&&+\mathbb{E}\left(\int_0^s\int_Z|\int_{\mathbb{R}^d}K(s,r,w)(g(r,x-w,z)-g(r,y-w,z))dw|^q\nu(dz)dr\right)\\
&\leq&N \mathbb{E}\left(\int_{\frac{3a-b}{2}}^s\int_Z\Big|\int_{\mathbb{R}^d}K(s,r,w)(g(r,x-w,z)-g(r,y-w,z))dw|^2\nu(dz)dr\right)^{q/2}\\
&&+N\mathbb{E}\left(\int_{\frac{3a-b}{2}}^s\int_Z|\int_{\mathbb{R}^d}K(s,r,w)(g(r,x-w,z)-g(r,y-w,z))dw|^q\nu(dz)dr\right)\\
&&+N \mathbb{E}\left(\int_0^{\frac{3a-b}{2}}\int_Z\Big|\int_{\mathbb{R}^d}(K(s,r,x-w)-K(s,r,y-w))g(r,w,z)dw|^2\nu(dz)dr\right)^{q/2}\\
&&+N\mathbb{E}\left(\int_0^{\frac{3a-b}{2}}\int_Z|\int_{\mathbb{R}^d}(K(s,r,x-w)-K(s,r,y-w))g(r,w,z)dw|^q\nu(dz)dr\right)\\
&\leq&N \mathbb{E}\left(\int_{\frac{3a-b}{2}}^s\int_Z\Big|\int_{\mathbb{R}^d}|K(s,r,w)g(r,x-w,z)|dw|^2\nu(dz)dr\right)^{q/2}\\
&&+N \mathbb{E}\left(\int_{\frac{3a-b}{2}}^s\int_Z\Big|\int_{\mathbb{R}^d}|K(s,r,w)g(r,y-w,z)|dw|^2\nu(dz)dr\right)^{q/2}\\
&&+N\mathbb{E}\left(\int_{\frac{3a-b}{2}}^s\int_Z|\int_{\mathbb{R}^d}|K(s,r,w)g(r,x-w,z)|dw|^q\nu(dz)dr\right)\\
&&+N\mathbb{E}\left(\int_{\frac{3a-b}{2}}^s\int_Z|\int_{\mathbb{R}^d}|K(s,r,w)g(r,y-w,z)|dw|^q\nu(dz)dr\right)\\
&&+N \mathbb{E}\left(\int_0^{\frac{3a-b}{2}}\int_Z\Big|\int_{\mathbb{R}^d}(K(s,r,x-w)-K(s,r,y-w))g(r,w,z)dw|^2\nu(dz)dr\right)^{q/2}\\
&&+N\mathbb{E}\left(\int_0^{\frac{3a-b}{2}}\int_Z|\int_{\mathbb{R}^d}(K(s,r,x-w)-K(s,r,y-w))g(r,w,z)dw|^q\nu(dz)dr\right)\\
&=:&I_{21}+\cdots+I_{26}.
   \eess
Similar to $I_{11}+I_{12}$, the four terms $I_{21}+\cdots+I_{24}$ is less than or equal to
   \bess
&&N[\varphi(2(b-a)c^{-\gamma})]^{q/q_0}\left(\left\|\left(\int_Z\|g(\cdot,\cdot,z)\|^2_{L^\infty(\mathcal {O}_T)}\nu(dz)\right)^{q/2}\right\|_{L^{\frac{\kappa}{\kappa-q}}}
\right.\\
&&\qquad\qquad\qquad\qquad\qquad\left.
+\left\|\int_Z\|g(\cdot,\cdot,z)\|^q_{L^\infty(\mathcal {O}_T)}\nu(dz)\right\|_{L^{\frac{\kappa}{\kappa-q}}}\right).
  \eess
Using {\bf Assumption 2.3 (iii)} with $\lambda=\frac{3a-b}{2}$, we get
   \bess
I_{25}+I_{26}&\leq&
N \mathbb{E}\left(\int_0^{\frac{3a-b}{2}}\Big|\int_{\mathbb{R}^d}
|K(s,r,x-w)-K(s,r,y-w)|dw\Big|^2dr\right.\\
&&\left.\times\int_Z\|g(\cdot,\cdot,z)\|^2_{L^\infty(\mathcal {O}_T)}\nu(dz)\right)^{q/2}\\
&&+N\mathbb{E}\left(\int_0^{\frac{3a-b}{2}}\Big|\int_{\mathbb{R}^d}|K(s,r,x-w)-K(s,r,y-w)|dw\Big|^qdr\right.\\
&&\left.\times
\int_Z\|g(\cdot,\cdot,z)\|^q_{L^\infty(\mathcal {O}_T)}\nu(dz)\right)\\
&\leq&N\varphi(2^{1+1/\gamma}c(b-a)^{-1/\gamma})^{q/q_0}\left(\left\|\left(\int_Z\|g(\cdot,\cdot,z)\|^2_{L^\infty(\mathcal {O}_T)}\nu(dz)\right)^{q/2}\right\|_{L^{\frac{\kappa}{\kappa-q}}}
\right.\\
&&\left.
+\left\|\int_Z\|g(\cdot,\cdot,z)\|^q_{L^\infty(\mathcal {O}_T)}\nu(dz)\right\|_{L^{\frac{\kappa}{\kappa-q}}}\right).
  \eess
Combining the above discussion, (\ref{3.2}) is obtained. The proof of this lemma is complete. $\Box$

Now, we are ready to prove the main result. The proof is similar to that of Theorem of 2.4 in \cite{Kim}.

{\bf Proof of Theorem \ref{t2.1}.}  Let $q\in[2,q_0\wedge \kappa]$. It suffices to prove that for each
  \bess
Q=Q_c(t_0,x_0):=(t_0-c^\gamma,t_0+c^\gamma)\times B_c(x_0)\subset\mathcal {O}_T,\ \ c>0,t_0>0,
  \eess
we have
   \bes
&&\frac{1}{Q^2}\mathbb{E}\int_Q\int_Q|\mathcal {G}g(t,x)-\mathcal {G}g(s,y)|^qdtdxdsdy\nm\\
&\leq&N\left(\left\|\left(\int_Z
\|g(\cdot,\cdot,z)\|^2_{L^\infty(\mathcal {O}_T)}\nu(dz)\right)^{q/2}\right\|_{L^{\frac{\kappa}{\kappa-q}}}
\right.\nm\\
&&\qquad+\left\|\int_Z|g(\cdot,\cdot,z)|^{q_0}_{L^\infty(\mathcal {O}_T)}\nu(dz)\right\|_{L^{\tilde\kappa}(\Omega)}^{q/q_0}\nm\\
&&\qquad\left.
+\left\|\int_Z\|g(\cdot,\cdot,z)\|^q_{L^\infty(\mathcal {O}_T)}\nu(dz)\right\|_{L^{\frac{\kappa}{\kappa-q}}}\right),
  \lbl{3.3}\ees
where $N=N(T,q,\varphi)$. Since the operator $\mathcal {G}$ is translation invariant with respect to $x$, i.e.
   \bess
\mathcal {G}g(\cdot,\cdot)(t,x+x_0)=\mathcal {G}g(\cdot,x_0+\cdot)(t,x),
   \eess
we may assume that $x_0=0$. We divide the left hand side of (\ref{3.3}) into two parts. Indeed,
   \bess
&&\frac{1}{Q^2}\mathbb{E}\int_Q\int_Q|\mathcal {G}g(t,x)-\mathcal {G}g(s,y)|^qdtdxdsdy\nm\\
&\leq&\frac{2}{Q}\mathbb{E}\int_Q|\mathcal {G}g_1(t,x)|^qdtdxdsdy\nm\\
&&+\frac{1}{Q^2}\mathbb{E}\int_Q\int_Q|\mathcal {G}g_2(t,x)-\mathcal {G}g_2(s,y)|^qdtdxdsdy\nm\\
&=:&J_1+J_2,
   \eess
where
  \bess
g_1(t,x,z)=I_{((t_0-2c^\gamma)\vee0,t_0+2c^\gamma)\times B_{2c}\times Z}(t,x,z)g(t,x,z), \ \ g_2=g-g_1.
  \eess

{\bf Estimate of $J_1$.} Since $Q\subset\mathcal {O}_T$, it holds that $t_0-c^\gamma\geq0$ and thus
   \bess
(t_0-c^\gamma,t_0+c^\gamma)\subset(t_0-2c^\gamma)\vee0,t_0+2c^\gamma)
   \eess
and $g$ vanishes on
   \bess
\Big[((t_0-2c^\gamma)\vee0,t_0+2c^\gamma)\times B^c_{2c}\times Z\Big]\bigcup\Big[(0,(t_0-2c^\gamma)\vee0)\times\mathbb{R}^d\times Z\Big].
  \eess
It follows from Lemma \ref{l3.1} with $a=(t_0-2c^\gamma)\vee0$ and $b=t_0+2c^\gamma$ that
   \bes
J_1\leq N\left\|\int_Z|g(\cdot,\cdot,z)|^{q_0}_{L^\infty(\mathcal {O}_T)}\nu(dz)\right\|_{L^{\tilde\kappa}(\Omega)}^{q/q_0}.
  \lbl{3.4}\ees

{\bf Estimate of $J_2$.} If $t_0\leq2c^\gamma$, we apply Lemma \ref{l3.2} with $a=t_0-c^\gamma$ and $b=t_0+c^\gamma$.
In this case, one can easily check that $bc^{-\gamma}\leq3$ and
   \bess
g_2=0\ \ \ \ {\rm on}\ \ \Big[(0,t_0+2c^\gamma)\times B_{2c}\times Z\Big].
   \eess
(\ref{3.1}) of Lemma \ref{l3.2} yields that
   \bes
J_2&\leq& N \left(\left\|\left(\int_Z
\|g(\cdot,\cdot,z)\|^2_{L^\infty(\mathcal {O}_T)}\nu(dz)\right)^{q/2}\right\|_{L^{\frac{\kappa}{\kappa-q}}}
\right.\nm\\
&&\left.
+\left\|\int_Z\|g(\cdot,\cdot,z)\|^q_{L^\infty(\mathcal {O}_T)}\nu(dz)\right\|_{L^{\frac{\kappa}{\kappa-q}}}\right)
  \lbl{3.5}\ees

On the other hand, if $t_0>2c^\gamma$, we apply Lemma \ref{l3.3} with $a=t_0-c^\gamma$ and $b=t_0+c^\gamma$.
In this case, one can easily check that $3a>b$ and
   \bess
g_2=0\ \ \ \ {\rm on}\ \ \Big[(t_0-2c^\gamma,t_0+2c^\gamma)\times B_{2c}\times Z\Big].
   \eess
Moreover, by using the nondecreasing of $\varphi$, we have
   \bess
&&\sup_{t_0\in\mathbb{R}_+,c>0}\Phi(t_0-c^\gamma,t_0+c^\gamma,c)\\
&=&\sup_{t_0\in\mathbb{R}_+,c>0}\left\{[\varphi(2)]^{q/q_0}+
[\varphi((b-a)c^{-\gamma})]^{q/q_0}+[\varphi(2^{1+1/\gamma}c(b-a)^{-1/\gamma})]^{q/q_0}\right\}\Big|_{\left\{\begin{array}{lll}a=t_0-c^\gamma\\
b=t_0+c^\gamma\end{array}\right.}\\
&<&\infty.
  \eess
(\ref{3.2}) implies that
   \bes
J_2&\leq& N\left(\left\|\left(\int_Z\|g(\cdot,\cdot,z)\|^2_{L^\infty(\mathcal {O}_T)}\nu(dz)\right)^{q/2}\right\|_{L^{\frac{\kappa}{\kappa-q}}}
\right.\nm\\
&&\left.
+\left\|\int_Z\|g(\cdot,\cdot,z)\|^q_{L^\infty(\mathcal {O}_T)}\nu(dz)\right\|_{L^{\frac{\kappa}{\kappa-q}}}\right).
  \lbl{3.6}\ees
Combining (\ref{3.4}), (\ref{3.5}) and (\ref{3.6}), we obtain (\ref{3.3}). The proof of Theorem \ref{t2.1} is complete. $\Box$

\begin{remark}\lbl{r3.1} In this paper, we only consider the simply case. Actually, one
can use the similar method and Kunita's second inequality (see Page 268 in \cite{Ap}) to deal with the
following case
    \bess
\mathcal {G}\hat g(t,x)&=&\int_0^t\int_{\mathbb{R}^d}K(t,s,x-y)h(s,y)dydW(s)\nm\\
&&+\int_0^t\int_Z\int_{\mathbb{R}^d}K(t,s,x-y)g(s,y,z)dy\tilde N(dz,ds),
    \eess
where $W$ and $\tilde N$ is a Wiener process and a compensated Poisson measure, respectively. Also see
\cite{MPR} for this case.
\end{remark}

\section{Applications}
\setcounter{equation}{0}

In this section, applying Theorem \ref{t2.1}, we obtain the BMO estimate of the following stochastic
singular integral operator
    \bes
 \mathcal {G} g(t,x)=\sum_{k=1}^\infty\int_0^t\int_{\mathbb{R}^m}\int_{\mathbb{R}^d}K(t,s,x-y)g^k(s,y)dyz\tilde N_k(dz,ds),
    \lbl{4.1}\ees
where $K(t,s,x)=\nabla^{\beta}p(t,s,x)$ and $p(t,s,x)$ is the heat kernel of the equation
   \bess
\partial_tu=\Delta^{\frac{\alpha}{2}}u.
  \eess
The fractional derivative of spatial variable is understood in sense of Fourier transform.
It is easy to see that
    \bess
\sum_{k=1}^\infty\int_0^t\int_{\mathbb{R}^m}\int_{\mathbb{R}^d}K(t,s,x-y)g^k(s,y)dyz\tilde N_k(dz,ds)
   \eess
is the fundamental solution to the following equation
   \bes
du_t(x)=\Delta^{\frac{\alpha}{2}}u_t(x)dt+\sum_{k=1}^\infty\int_{\mathbb{R}^m}g^k(t,x)z\tilde N_k(dz,dt),\ \ \ u_0=0,\ 0\leq t\leq T,
  \lbl{0.1}\ees
where $\int_{\mathbb{R}^m}z\tilde N_k(t,dz)=:Y_t^k$ are independent $m$-dimensional pure jump L\'{e}vy processes with L\'{e}vy measure of $\nu^k$.
Indeed, one can use the method of \cite{Kim} (see the proof of Lemma 6.1) to prove the above result.
On the other hand, Kim-Kim \cite{KpK} considered the general case. We only recall the results
concerned with this paper. In section 3 of \cite{KpK}, Kim-Kim studied the following linear equation (see Page 3935 of \cite{KpK}):
   \bes
du=(a(\omega,t)\Delta^{\frac{\alpha}{2}}u+f)dt+\sum_{i=1}^\infty h^kdW_t^k+\sum_{k=1}^{\infty}\sum_{j=1}^mg^{k,j}\cdot dY_t^{k,j},\ \ u(0)=u_0,
   \lbl{4.2}\ees
where $h=(h^1,h^2,\cdots)$, $W_t^k$ is independent one-dimensional Wiener processes and $Y_t^k:=\int_{\mathbb{R}^m}z\tilde N_k(t,dz)$. Note
that $Y_t^k$ are independent $m$-dimensional pure jump L\'{e}vy processes with L\'{e}vy measure of $\nu^k$.
For any $q,k=1,2,\cdots$, denote
   \bess
\hat c_{k,q}:=\left(\int_{\mathbb{R}^m}|z|^q \nu^k(dz)\right)^{\frac{1}{q}}.
   \eess
Fix $p\in[2,\infty)$ and set $\hat c_k:=\hat c_{k,2}\vee\hat c_{k,p}$. Assume that
   \bes
\hat c:=\sup_{k\geq1}\hat c_k<\infty.
   \lbl{4.3}\ees

Let $\mathcal {P}$ be the predictable $\sigma$-field generated by $\{\mathcal {F}_t,t\geq0\}$ and
$\bar{\mathcal {P}}$ be the completion of $\mathcal {P}$ with respect to $dP\times dt$. For $\eta\in\mathbb{R}$, define
$\mathbb{H}^\eta_p(T):=L^p(\Omega\times[0,T],\bar{\mathcal {P}},H_p^\eta)$, that is, $\mathbb{H}^\eta_p(T)$ is the set of
all $\bar{\mathcal {P}}$-measurable processes $u:\Omega\times[0,T]\mapsto H_p^\eta$ so that
   \bess
\|u\|_{\mathbb{H}^\eta_p(T)}:=\left(\mathbb{E}\int_0^T\||u(\omega,t,\cdot)\|^p_{H_p^\eta}dt\right)^{1/p}<\infty,
   \eess
where $H_p^\eta(\mathbb{R}^d):=\{u:\,D^{\textbf{n}}u\in L^p(\mathbb{R}^d),\,|\textbf{n}|\leq \eta\}$ for $\eta=1,2,\dots$.
And when $\eta$ is not an integer, $H_p^\eta(\mathbb{R}^d)$ is defined by Fourier transform.

For $\ell_2$-valued $\bar{\mathcal {P}}$-measurable processes $g=(g^1,g^2,\cdots)$, we write $g\in\mathbb{H}^\eta_p(T,\ell_2)$ if
    \bess
\|g\|_{\mathbb{H}^\eta_p(T,\ell_2)}&:=&\left(\mathbb{E}\int_0^T\|g(\omega,t,\cdot)\|^p_{H^\eta_p(T,\ell_2)}dt\right)^{1/p}\\
&=&\left(\mathbb{E}\int_0^T\||(1-\Delta)^{\eta/2}g(\omega,t,\cdot)|_{\ell_2}\|^p_pdt\right)^{1/p}<\infty.
   \eess

Lastly, we define
   \bess
\|u\|_{\mathcal {H}_p^{\eta+\alpha}(T)}&:=&\|u\|_{\mathbb{H}^{\eta+\alpha}_p(T)}+\|f\|_{\mathbb{H}^{\eta+\alpha}_p(T)}
+\|h\|_{\mathbb{H}^{\eta+\alpha/2}_p(T,\ell_2)}\\
&&+\sum_{j=1}^m\|g^{\cdot,j}\|_{\mathbb{H}^{\eta+\alpha/2}_p(T,\ell_2)}
+\|u(0)\|_{U_p^{\eta+\alpha-\alpha/p}},
   \eess
where $\|u(0)\|_{U_p^{\eta+\alpha-\alpha/p}}:=\left(\mathbb{E}[\|u_0\|_{H^\eta_p}^p]\right)^{1/p}$.

\begin{prop}\lbl{p4.1} {\rm\cite[Theorem 3.6]{KpK}} Suppose (\ref{4.3}) holds. Then for any $f\in \mathcal {H}_p^{\eta}(T)$,
$h\in \mathbb{H}^{\eta+\alpha/2}_p(T,\ell_2)$, $g^{\cdot,j}\in\mathbb{H}^{\eta+\alpha-\alpha/p}_p(T,\ell_2) $,
$1\leq j\leq m$ and $u_0\in U_p^{\eta+\alpha-\alpha/p}$, Eq. (\ref{4.2}) has a unique solution $u$ in $\mathcal {H}_p^{\eta+\alpha}$,
and for this solution
   \bess
\|u\|_{\mathcal {H}_p^{\eta+\alpha}(t)}&\leq &N(p,T,a)\Big(\|f\|_{\mathbb{H}^{\eta}_p(t)}+\|h\|_{\mathbb{H}^{\eta+\alpha/2}_p(t,\ell_2)}\\
&&+\sum_{j=1}^m\|g^{\cdot,j}\|_{\mathbb{H}^{\eta+\alpha-\alpha/p}_p(t,\ell_2)}
+\|u(0)\|_{U_p^{\eta+\alpha-\alpha/p}}\Big)
   \eess
for every $t\leq T$.
\end{prop}

In order to investigate the BMO estimate of the solution, we recall some properties of kernel $p(t,s,x)$
(see \cite{B-J, B-S-S, CH, I} for more details).
 \begin{itemize}
 \item for any $t>0$,
\bess
\|p(t, \cdot)\|_{L^1(\R^d)}=1 \text{ for all } t>0.
\eess
 \item $ p(t, x, y)$ is  $C^\infty$ on $(0,\infty)\times \R^d\times \R^d$ for each $t>0$;

 \item   for $t>0$, $x, y\in \R^d$, $x\neq y$, the sharp estimate of $\widehat p(t, x)$ is
\bess
p(t, x, y)\approx \min\left( \frac{t}{|x-y|^{d+\alpha}}, t^{-d/\alpha}\right);
  \eess
\item  for $t>0$, $x, y\in \R^d$, $x\neq y$,   the  estimate of the first order derivative of $\widehat p(t, x)$ is
\bes
  |\nabla_x p(t, x, y)|\approx  |y-x|\min\left\{ \frac{t}{|y-x|^{d+2+\alpha}}, t^{-\frac{d+2}{\alpha}}\right\}.
 \lbl{4.4}\ees
\end{itemize}
The notation $f(x)\approx g(x)$ means that there is a number $0<C<\infty$ independent of $x$,
i.e. a constant, such that for every $x$ we have $C^{-1}f(x)\leq g(x)\leq Cf(x)$.
The estimate (\ref{4.4}) for the first order derivative of $ p(t,x)$ was derived in  \cite[Lemma 5]{B-J}.
Xie et al. \cite{XDLL} the estimate of the $m$-th order derivative of $p(t, x)$ by induction.
\begin{prop}
\label{p4.2}{\rm\cite[Lemma 2.1]{XDLL}}
For any $m\geq 0$, we have
  \bes
\partial_x^m p(t, x)=\sum_{n=0}^{n=\lfloor \frac{m}{2}\rfloor}C_n |x|^{m-2n} \min \left\{ \frac{t}{|x|^{d+\alpha+2(m-n)}}, t^{-\frac{d+2(m-n)}{\alpha}}\right\},
    \lbl{4.5}\ees
where $\lfloor \frac m 2\rfloor$ means the largest integer that is less  than $\frac m2$.
\end{prop}

Next, we claim that the kernel $\nabla^{\frac{\alpha}{q_0}}p(t,s,x)$, $q_0\geq2$, satisfies the {\bf Assumption 2.3} with $\gamma=\alpha$
and $\kappa=\infty$.

\begin{lem}\lbl{l4.1}Let $\beta=\frac{\alpha}{q_0}$. The following estimates hold.

 (i) For any $t>\lambda>0$ and $c>0$,
   \bess
\int_\lambda^t\Big|\int_{|x|\geq c}|\nabla^{\beta}p(t,r,x)|dx\Big|^{q_0}dr\leq N\left( [(t-\lambda)c^{-\alpha}]^{q_0+1}+[(t-\lambda)c^{-\alpha}]\right);
   \eess

(ii) For any $t>s>\lambda>0$,
   \bess
\int_0^\lambda\left(\int_{\mathbb{R}^d}|\nabla^\beta p(t,r,x)-\nabla^{\beta}p(s,r,x)|dx\right)^{q_0}dr
\leq  N[(t-s)(t\wedge s-\lambda)^{-1}]^{q_0};
   \eess

(iii) For any $s>\lambda\geq0$ and $h\in\mathbb{R}^d$,
    \bess
\int_0^\lambda\left(\int_{\mathbb{R}^d}|\nabla^{\beta}p(s,r,x+h)-\nabla^{\beta}p(s,r,x)|dx\right)^{q_0}dr
\leq N\varphi(|h|(s-\lambda)^{-1/\alpha}).
   \eess
\end{lem}

{\bf Proof.} Note that $\beta=\frac{\alpha}{q_0}<2$. By using Proposition \ref{p4.2}, we have if $c>(t-r)^{\frac{1}{\alpha}}$,
   \bess
&&\int_\lambda^t\Big|\int_{|x|\geq c}|\nabla^{\beta}p(t,r,x)|dx\Big|^{q_0}dr\\
&\leq&N\int_\lambda^t\Big|\int_{|x|\geq c}|x|^{\beta}\frac{t-r}{|x|^{d+\alpha+2\beta}}dx\Big|^{q_0}dr\\
&\leq&N\int_\lambda^t\Big|\int_c^\infty|x|^{\beta}\cdot|x|^{d-1}\frac{t-r}{|x|^{d+\alpha+2\beta}}d|x|\Big|^{q_0}dr\\
&=&Nc^{-\alpha(q_0+1)}\int_\lambda^t(t-r)^{q_0}dr\\
&\leq&N [(t-\lambda)c^{-\alpha}]^{q_0+1}.
   \eess
When $c\leq(t-r)^{\frac{1}{\alpha}}$, we have $(t-r)^{-1}\leq c^{-\alpha}$
   \bess
&&\int_\lambda^t\Big|\int_{|x|\geq c}|\nabla^{\beta}p(t,r,x)|dx\Big|^{q_0}dr\\
&\leq&N\int_\lambda^t\left(\int_{(t-r)^{\frac{1}{\alpha}}}^\infty|x|^{\beta}\cdot|x|^{d-1}\frac{t-r}{|x|^{d+\alpha+2\beta}}d|x|\right.\\
&&\qquad\left.+\int_c^{(t-r)^{\frac{1}{\alpha}}}|x|^{\beta}\cdot|x|^{d-1}(t-r)^{-\frac{d+2\beta}{\alpha}}d|x|
\right)^{q_0}dr\\
&\leq&N\int_\lambda^t\left(\int_c^\infty|x|^{\beta}\cdot|x|^{d-1}\frac{t-r}{|x|^{d+\alpha+2\beta}}d|x|\right.\\
&&\qquad\left.+\int_0^{(t-r)^{\frac{1}{\alpha}}}|x|^{\beta}\cdot|x|^{d-1}(t-r)^{-\frac{d+2\beta}{\alpha}}d|x|
\right)^{q_0}dr\\
&\leq&Nc^{-\alpha(q_0+1)}\int_\lambda^t(t-r)^{q_0}dr+Nc^{-\alpha}\int_\lambda^tdr\\
&\leq&N [(t-\lambda)c^{-\alpha}]^{q_0+1}+[(t-\lambda)c^{-\alpha}].
   \eess
Hence we obtain the first estimate.

When $\alpha+\frac{\alpha}{q_0}<2$, $\lfloor \frac{\alpha+\alpha/q_0}{2}\rfloor=0$.
Using the fact that $\partial_tp=\Delta^{\alpha/2}p$, $\beta q_0=1$ and
Proposition \ref{p4.2}, we get
   \bess
&&\int_0^\lambda\left(\int_{\mathbb{R}^d}|\nabla^{\beta}p(t,r,x)-\nabla^{\beta}p(s,r,x)|dx\right)^{q_0}dr\\
&\leq&(t-s)^{q_0}\int_0^\lambda\left(\int_{\mathbb{R}^d}|\nabla^{\alpha+\beta}p(\xi-r,x)|dx\right)^{q_0}dr\\
&\leq&N (t-s)^{q_0}\int_0^\lambda\left(\int_0^{(\xi-r)^{\frac{1}{\alpha}}}|x|^{\alpha+\beta}|x|^{d-1}(\xi-r)^{-\frac{d+2\alpha+2\beta}{\alpha}}d|x|\right.\\
&&\qquad\qquad\qquad\left.+\int_{(\xi-r)^{\frac{1}{\alpha}}}^\infty|x|^{\alpha+\beta}|x|^{d-1}\frac{\xi-r}{|x|^{d+3\alpha+2\beta}}d|x|\right)^{q_0}dr\\
&\leq&N(t-s)^{q_0}\int_0^\lambda(\xi-r)^{-q_0-1}dr\\
&\leq& N[(t-s)(t\wedge s-\lambda)^{-1}]^{q_0},
  \eess
where $\xi=\theta t+(1-\theta)s$, $\theta\in[0,1]$.

Since $q_0\geq2$ and $0\leq\alpha\leq2$, we have $\alpha+\frac{\alpha}{q_0}<4$. When $2\leq\alpha+\frac{\alpha}{q_0}<4$, we have
   \bess
&&\int_0^\lambda\left(\int_{\mathbb{R}^d}|\nabla^{\beta}p(t,r,x)-\nabla^{\frac{\alpha}{q_0}}p(s,r,x)|dx\right)^{q_0}dr\\
&\leq&(t-s)^{q_0}\int_0^\lambda\left(\int_{\mathbb{R}^d}|\nabla^{\alpha+\beta}p(\xi-r,x)|dx\right)^{q_0}dr\\
&\leq&N (t-s)^{q_0}\int_0^\lambda\left(\int_0^{(\xi-r)^{\frac{1}{\alpha}}}|x|^{\alpha+\beta}|x|^{d-1}(\xi-r)^{-\frac{d+2\alpha+2\beta}{\alpha}}d|x|\right.\\
&&\qquad\qquad\qquad\int_0^{(\xi-r)^{\frac{1}{\alpha}}}|x|^{\alpha+\beta-2}|x|^{d-1}(\xi-r)^{-\frac{d+2\alpha+2\beta-2}{\alpha}}d|x|\\
&&\qquad\qquad\qquad+\int_{(\xi-r)^{\frac{1}{\alpha}}}^\infty|x|^{\alpha+\beta}|x|^{d-1}\frac{\xi-r}{|x|^{d+3\alpha+2\beta}}d|x|\\
&&\qquad\qquad\qquad\left.+\int_{(\xi-r)^{\frac{1}{\alpha}}}^\infty|x|^{\alpha+\beta-2}|x|^{d-1}\frac{\xi-r}{|x|^{d+3\alpha+2\beta-2}}d|x|\right)^{q_0}dr\\
&\leq&N(t-s)^{q_0}\int_0^\lambda(\xi-r)^{-q_0-1}dr\\
&\leq& N[(t-s)(t\wedge s-\lambda)^{-1}]^{q_0},
  \eess
where $\xi=\theta t+(1-\theta)s$, $\theta\in[0,1]$. Thus we obtain the second estimate.

For the last estimate (iii), noting that $1+\beta\leq2$, we have for $1+\beta<2$
   \bess
&&\int_0^\lambda\left(\int_{\mathbb{R}^d}|\nabla^{\beta}p(s,r,x+h)-\nabla^{\beta}p(s,r,x)|dx\right)^{q_0}dr\\
&\leq&N \int_0^\lambda h^{q_0}\left(\int_{\mathbb{R}^d}|\nabla^{1+\beta}p(s,r,x+\theta h)|dx\right)^{q_0}dr\\
&\leq&N \int_0^\lambda h^{q_0}\left(\int_0^{(s-r)^{\frac{1}{\alpha}}}|x|^{1+\beta}\cdot|x|^{d-1}(s-r)^{-\frac{d+2+2\beta}{\alpha}}d|x|\right.\\
&&\qquad\qquad\qquad\left.+\int_{(s-r)^{\frac{1}{\alpha}}}^\infty|x|^{1+\beta}\cdot|x|^{d-1}\frac{s-r}{|x|^{d+\alpha+2+2\beta}}d|x|\right)^{q_0}dr\\
&\leq&N[h(s-\lambda)^{-1}]^{q_0},
  \eess
where $\theta\in[0,1]$. When $1+\beta=2$, similar the case (ii), one can get the same estimate.
The proof of Lemma is complete. $\Box$

It follows from the Proposition \ref{p4.1} that $\nabla^\beta p(t,s,x)$ satisfies the {\bf Assumption 2.4}.
By using Theorem \ref{t2.1}, we have the following result.
\begin{theo}\lbl{t4.1} Let $q_0\geq2$. Suppose (\ref{4.3}) with $p\geq q_0$ holds. Then for any $g\in\mathbb{H}^{\eta+\alpha-\alpha/p}_p(T,\ell_2) $,
Eq. (\ref{0.1}) has a unique solution $u$ in $\mathcal {H}_p^{\eta+\alpha}$ ($\eta\in\mathbb{R}$),
and for this solution
   \bess
\|u\|_{\mathcal {H}_p^{\eta+\alpha}(t)}&\leq &N(p,T)\|g\|_{\mathbb{H}^{\eta+\alpha-\alpha/p}_p(t,\ell_2)}
   \eess
for every $t\leq T$.

Moreover, we have for $q\in[2,q_0]$
      \bess
[\nabla^{\beta}u]_{\mathbb{BMO}(T,q)}&\leq& N\hat c\left(\mathbb{E}[\||g|_{\ell_2}\|^{q_0}_{L^\infty(\mathcal {O}_T)}]\right)^{q/q_0},
  \eess
where $\beta=\alpha/q_0$ and $\hat c$ is  defined as in (\ref{4.3}).
\end{theo}

When the L\'{e}vy noise is replaced by Brownian motion in (\ref{0.1}), i.e.,
   \bes
du_t(x)=\Delta^{\frac{\alpha}{2}}u_t(x)dt+\sum_{k=1}^\infty h^k(t,x)dW_t^k,\ \ \ u_0=0,\ 0\leq t\leq T,
  \lbl{4.7}\ees
where $W_t^k$ are independent one-dimensional Wiener processes. Denote $h=(h^1,h^2,\cdots)$.

Similar to Lemma \ref{l4.1}, one can prove $\nabla^{\frac{\alpha}{2}}p(t,s,x)$
satisfies the {\bf Assumption 2.1}. From Proposition \ref{p4.1}, we know that
{\bf Assumption 2.2} holds for $\nabla^{\frac{\alpha}{2}}p(t,s,x)$.
Thus we can get the following result.
\begin{theo}\lbl{t4.2}  Suppose that $h\in L^p(T,\ell_2)$, there exists a uniqueness
solution $u$ in $\mathcal {H}_p^{\eta+\alpha}$ ($\eta\in\mathbb{R}$),
and for this solution
   \bess
\|u\|_{\mathcal {H}_p^{\eta+\alpha}(t)}&\leq &N(p,T)\|h\|_{\mathbb{H}^{\eta+\alpha/2}_p(t,\ell_2)}
   \eess
for every $t\leq T$.

Moreover, we have for any $q\in(0,p]$
      \bess
[\nabla^{\frac{\alpha}{2}}u]_{\mathbb{BMO}(T,q)}&\leq& N\left(\mathbb{E}[\||h|_{\ell_2}\|^{p}_{L^\infty(\mathcal {O}_T)}]\right)^{1/p}.
  \eess
\end{theo}

\begin{remark}\lbl{r4.1} 1. In Lemma \ref{l4.1}, the second part (ii) is essential.
From the proof of Theorem \ref{t2.1}, the bound of the BMO norm can be controlled
by the function $\varphi$ and some norm of $g$, where the bound of the function $\varphi$
depends on the choice of scale of time and space. In second part (ii),
we must prove that the left hand side of (ii) can be
controlled by  the function of $(t-s)(t\wedge s-\lambda)^{-1}$. Only in this form,
the left hand side of (ii) can be
controlled by a constant.

2. Particularly, taking $q_0=2$, we have Lemma \ref{l4.1} holds for $\nabla^{\frac{\alpha}{2}}p(t,s,x)$.
Hence we have Theorem \ref{t4.2}. Noting that if $\alpha=2$, Theorem \ref{t4.2} becomes
\cite[Theorem 3.4]{Kim}. Thus we generalize the result of \cite{Kim}.
  \end{remark}

\section{Discussion}
\setcounter{equation}{0}

In this section, we give another proof of Theorem \ref{t2.1} under some assumptions on $g$.
Similarly, one can give another proof of \cite[ Theorem 2.4]{Kim} under the same assumptions on $g$.
Firstly, let us recall the proofs of Theorem \ref{t2.1} and \cite[ Theorem 2.4]{Kim}.
The reason why we divide
the interval $(0,s)$ into two parts
$(0,\frac{3a-b}{2})$ and $(\frac{3a-b}{2},s)$ in proof of Lemma \ref{3.3} is the singularity of $K$ at time $t$.
In order to see it clearly, we look at the Section 4 and recall that for any $t>\lambda>0$ and $c>0$,
   \bess
\int_\lambda^t\Big|\int_{|x|\geq c}|\nabla^{\beta}p(t,r,x)|dx\Big|^{q_0}dr\leq N\left( [(t-\lambda)c^{-\alpha}]^{q_0+1}+[(t-\lambda)c^{-\alpha}]\right).
   \eess
Note that if we choose $c=0$, then the above integral will be infinity. Indeed, direct calculations
show that
  \bess
 \int_\lambda^t\Big|\int_{\mathbb{R}^d}|\nabla^{\beta}p(t,r,x)|dx\Big|^{q_0}dr
 \approx N\int_\lambda^t(t-r)^{-1}dr=\infty.
   \eess
Obviously, the singularity of $\nabla^{\beta}p$ appears at $t$. But $p\in L^1(\mathbb{R}^d)$,
thus a natural question appears: when the singularity of $p$ does not appear at $t$, is there another proof ?
Moreover, it is easy to see that the derivative of $p$ deduces the singularity of $\nabla^{\beta}p$ at $t$.
In this section, we first give a similar theorem to Theorem \ref{t2.1} under different assumptions.
Then as an application, we
use the method of integration by part to deal with the derivative of $p$ and obtain the BMO estimate
by direct calculation.

\begin{theo}\lbl{t5.1} Assume that the kernel function is a deterministic function and satisfies that for all $t\geq r\geq0$,
   \bess
\int_0^t\int_{\mathbb{R}^d}|K(t,r,x)|dxdr\leq N(T).
   \eess
Assume further that there exists a positive constant $q_0>2$ such that
   \bess
\mathbb{E}\left(\int_Z
\|g(\cdot,\cdot,z)\|^\varpi_{L^\infty(\mathcal {O}_T)}\nu(dz)\right)^{\frac{q_0}{2}}<\infty,\ \ \ \varpi=2\ {\rm or}\ q_0.
  \eess
Then for any $q\in(0,q_0]$,
one has
     \bess
[\mathcal {G}g]_{\mathbb{BMO}(T,q)}&\leq&
N\mathbb{E}\left(\int_Z
\|g(\cdot,\cdot,z)\|^2_{L^\infty(\mathcal {O}_T)}\nu(dz)\right)^{\frac{q}{2}}\\
&&+\mathbb{E}\left(\int_Z
\|g(\cdot,\cdot,z)\|^{q}_{L^\infty(\mathcal {O}_T)}\nu(dz)\right),
   \eess
where $N=N(N_0,d,q,q_0,T)$.
\end{theo}

{\bf Proof.} It suffices to prove that for each
  \bess
Q=Q_c(t_0,x_0):=(t_0-c^\gamma,t_0+c^\gamma)\times B_c(x_0)\subset\mathcal {O}_T,\ \ c>0,t_0>0,
  \eess
we have
   \bes
&&\frac{1}{Q^2}\mathbb{E}\int_Q\int_Q|\mathcal {G}g(t,x)-\mathcal {G}g(s,y)|^qdtdxdsdy\nm\\
&\leq&N\mathbb{E}\left(\int_Z
\|g(\cdot,\cdot,z)\|^2_{L^\infty(\mathcal {O}_T)}\nu(dz)\right)^{\frac{q}{2}}+\mathbb{E}\left(\int_Z
\|g(\cdot,\cdot,z)\|^{q}_{L^\infty(\mathcal {O}_T)}\nu(dz)\right),
  \lbl{5.1}\ees
where $N=N(T,q,\varphi)$. Since the operator $\mathcal {G}$ is translation invariant with respect to $x$,
we may assume that $x_0=0$. Kunita's first inequality implies that
    \bess
\mathbb{E}|\mathcal {G}g(t,x)|^q&\leq&
\mathbb{E}\left(\int_0^t\int_Z|\int_{\mathbb{R}^d}k(t-r,y)g(r,x-y,z)dy|^2\nu(dz)dr\right)^{q/2}\\
&&+\mathbb{E}\left(\int_0^t\int_Z|\int_{\mathbb{R}^d}k(t-r,y)g(r,x-y,z)dy|^q\nu(dz)dr\right)\\
&\leq&\mathbb{E}\left(\int_Z
\|g(\cdot,\cdot,z)\|^2_{L^\infty(\mathcal {O}_T)}\nu(dz)\times\int_0^t|\int_{\mathbb{R}^d}k(t-r,y)dy|^2dr\right)^{\frac{q}{2}}\\
&&+\mathbb{E}\left(\int_Z
\|g(\cdot,\cdot,z)\|^{q}_{L^\infty(\mathcal {O}_T)}\nu(dz)\times\int_0^t|\int_{\mathbb{R}^d}k(t-r,y)dy|^qdr\right)\\
&\leq&N(T)\mathbb{E}\left(\int_Z
\|g(\cdot,\cdot,z)\|^2_{L^\infty(\mathcal {O}_T)}\nu(dz)\right)^{\frac{q}{2}}\\
&&+\mathbb{E}\left(\int_Z
\|g(\cdot,\cdot,z)\|^{q}_{L^\infty(\mathcal {O}_T)}\nu(dz)\right)\\
&<&\infty.
  \eess
Thus we have
   \bess
&&\frac{1}{Q^2}\mathbb{E}\int_Q\int_Q|\mathcal {G}g(t,x)-\mathcal {G}g(s,y)|^qdtdxdsdy\nm\\
&\leq&\frac{2}{Q}\mathbb{E}\int_Q|\mathcal {G}g(t,x)|^qdtdx\\
&\leq&N(T)\mathbb{E}\left(\int_Z
\|g(\cdot,\cdot,z)\|^2_{L^\infty(\mathcal {O}_T)}\nu(dz)\right)^{\frac{q}{2}}\\
&&+\mathbb{E}\left(\int_Z
\|g(\cdot,\cdot,z)\|^{q}_{L^\infty(\mathcal {O}_T)}\nu(dz)\right),
   \eess
which implies (\ref{5.1}) holds.
The proof of Theorem \ref{t5.1} is complete. $\Box$

As an application, for simplicity, we consider the following stochastic evolution equation
    \bes
du=\Delta udt+\int_Zg(t,x,z)\tilde N(dt,dz)\ \ \ u(0,x)=0.
  \lbl{5.2}\ees
It is easy to check that the solution of (\ref{5.2}) is
   \bess
u(t,x)=\int_0^t\int_Z\int_{\mathbb{R}^d}K(t-r,y)g(r,y,z)dyd\tilde N(dr,dz).
   \eess
It follows the properties of heat kernel that
   \bess
\int_{\mathbb{R}^d}|K(t,r,x)|dx=1 \ \ {\rm for\ all }\ t>r>0.
   \eess
Applying Theorem \ref{t5.1}, we have
\begin{theo}\lbl{t5.2}
Assume that there exists a positive constant $q_0>2$ such that
   \bess
\mathbb{E}\left(\int_Z
\|g(\cdot,\cdot,z)\|^\varpi_{L^\infty(\mathcal {O}_T)}\nu(dz)\right)^{\frac{q_0}{2}}<\infty,\ \ \ \varpi=2\ {\rm or}\ q_0.
  \eess
Then for any $q\in(0,q_0]$,
one has
     \bess
[u]_{\mathbb{BMO}(T,q)}&\leq&
N\mathbb{E}\left(\int_Z
\|g(\cdot,\cdot,z)\|^2_{L^\infty(\mathcal {O}_T)}\nu(dz)\right)^{\frac{q}{2}}\\
&&+\mathbb{E}\left(\int_Z
\|g(\cdot,\cdot,z)\|^{q}_{L^\infty(\mathcal {O}_T)}\nu(dz)\right),
   \eess
where $N=N(N_0,d,q,q_0,T)$. Moreover, if we further assume that
   \bess
\mathbb{E}\left(\int_Z
\|\nabla_xg(\cdot,\cdot,z)\|^\varpi_{L^\infty(\mathcal {O}_T)}\nu(dz)\right)^{\frac{q_0}{2}}<\infty,\ \ \ \varpi=2\ {\rm or}\ q_0.
  \eess
Then for any $q\in(0,q_0]$,
one has
     \bess
[\nabla u]_{\mathbb{BMO}(T,q)}&\leq&
N\mathbb{E}\left(\int_Z
\|\nabla_xg(\cdot,\cdot,z)\|^2_{L^\infty(\mathcal {O}_T)}\nu(dz)\right)^{\frac{q}{2}}\\
&&+\mathbb{E}\left(\int_Z
\|\nabla_xg(\cdot,\cdot,z)\|^{q}_{L^\infty(\mathcal {O}_T)}\nu(dz)\right),
   \eess
where $N=N(N_0,d,q,q_0,T)$ and $\nabla_xg=\nabla_xg(t,\cdot,z)$.
\end{theo}

{\bf Proof.} Denote $u(t,x)=\mathcal {G}g(t,x)$. Noting that
   \bess
\nabla_x\mathcal {G}g(t,x)=\int_0^t\int_Z\int_{\mathbb{R}^d}k(t-r,y)\nabla_xg(r,x-y,z)dy\tilde N(dr,dz).
   \eess
Then similar to the proof of Theorem \ref{t5.1}, one can get the desired result. $\Box$

\begin{remark}\lbl{r5.1}  Comparing with the proofs of Theorems \ref{t2.1} and \ref{t5.1}, we find if
we assume the function $g$ has high regularity, then the proof of BMO estimate will be very
simple.
The proof of Theorem 4.1 will be also simple if we improve the regularity of $g$.

If $g\equiv0$, then $u\equiv0$. That is to say, the noise has effect on the
regularity of the solutions.
\end{remark}

\medskip

\noindent {\bf Acknowledgment} The first author was supported in part
by NSFC of China grants 11301146, 11171064.

 \end{document}